\documentstyle{amltd2004}
\begin{document}
\annalsline{155}{2002}
\startingpage{599}
\def\bye{\end{document}}
 \font\tenrm=cmr10
\def\joinrel{\mathrel{\mkern-4mu}}
\def\relbar{\mathrel{\smash-}}
\def\lrar{\relbar\joinrel\relbar\joinrel\relbar\joinrel\relbar\joinrel\rightarrow}
\def\llrar{\relbar\joinrel\relbar\joinrel\relbar\joinrel\relbar\joinrel\relbar\joinrel\rightarrow}
\newcommand{\Mot}{{\bf Mot}}
\newcommand{\Pic}{\mathop{\rm Pic}\nolimits}
\newcommand{\Hom}{\mathop{\rm Hom}\nolimits}
\newcommand{\End}{\mathop{\rm End}\nolimits}
\newcommand{\Ind}{\mathop{\rm Ind}\nolimits}
\newcommand{\Tr}{\mathop{\rm Tr}\nolimits}
\newcommand{\ob}{\mathop{\rm ob}\nolimits}
\newcommand{\ch}{\mathop{\rm char}\nolimits}
\newcommand{\Ker}{\mathop{\rm Ker}\nolimits}
\def\LCM{{\bf LCM}}
\def\LMot{{\bf LMot}}
%\input BoxedEps.Tex %for mac
%\SetTexturesEPSFSpecial %for mac
\input boxedeps.tex % unix
\SetepsfEPSFSpecial % unix
\HideDisplacementBoxes
\def\figin#1#2{\medbreak
$$
 {\BoxedEPSF{#1.eps scaled
#2}%
}%
$$
\medbreak\noindent}
\def\ev{{\rm ev}}
\newcommand{\Aut}{\mathop{\rm Aut}\nolimits}
\def\1{{1\mkern-7mu1}}
%--------------- Author macros ---------------
%for Bbb in amstex
\catcode`\@=11
\font\twelvemsb=msbm10 scaled 1100
\font\tenmsb=msbm10
%\font\ninemsb=msbm7 scaled 1100%msbm9
\font\ninemsb=msbm10 scaled 800
\newfam\msbfam
\textfont\msbfam=\twelvemsb  \scriptfont\msbfam=\ninemsb
  \scriptscriptfont\msbfam=\ninemsb
\def\msb@{\hexnumber@\msbfam}
\def\Bbb{\relax\ifmmode\let\next\Bbb@\else
 \def\next{\errmessage{Use \string\Bbb\space only in math
mode}}\fi\next}
\def\Bbb@#1{{\Bbb@@{#1}}}
\def\Bbb@@#1{\fam\msbfam#1}
\catcode`\@=12

 \catcode`\@=11
\font\twelveeuf=eufm10 scaled 1100
\font\teneuf=eufm10
\font\nineeuf=eufm7 scaled 1100%eufm9
\newfam\euffam
\textfont\euffam=\twelveeuf  \scriptfont\euffam=\teneuf
  \scriptscriptfont\euffam=\nineeuf
\def\euf@{\hexnumber@\euffam}
\def\frak{\relax\ifmmode\let\next\frak@\else
 \def\next{\errmessage{Use \string\frak\space only in math
mode}}\fi\next}
\def\frak@#1{{\frak@@{#1}}}
\def\frak@@#1{\fam\euffam#1}
\catcode`\@=12
%-------------- Author entries --------------------

\newcommand{\tfrac}{\frac}
\title{Polarizations and\\ Grothendieck's
standard conjectures} 
\shorttitle{Grothendieck's standard conjectures} 
  \acknowledgements{Part of this research was supported by the National Science
Foundation.}
\author{J. S. Milne}
\institutions{2679 Bedford Rd., Ann Arbor, MI 48104\\
{\eightpoint {\it E-mail address\/}: math@jmilne.org}\\
{\eightpoint {\it URL address\/}: www.jmilne.org/math/} }
%-------------- Article Text--------------------

 \centerline{\bf Abstract}
\vglue12pt
We prove that Grothendieck's Hodge standard conjecture holds for abelian
varieties in arbitrary characteristic if the Hodge conjecture holds for
complex abelian varieties of CM-type. For abelian varieties with no exotic
algebraic classes, we prove the Hodge standard conjecture unconditionally.

\vglue12pt\def\sni#1{\smallbreak\noindent{#1}. }
\def\ssni#1{\vglue-1pt\noindent\hskip18pt {#1}.}

\sni{1} Polarizations on quotient categories
\sni{2} Polarizations on categories of motives over finite fields
\sni{3} The Hodge standard conjecture
\smallbreak\noindent References

\intro

In examining Weil's proofs (Weil 1948) of the Riemann hypothesis for curves
and abelian varieties over finite fields, Grothendieck was led to state two
``standard'' conjectures (Grothendieck 1969), which imply the Riemann
hypothesis for all smooth projective varieties over a finite field,
essentially by Weil's original argument. Despite Deligne's proof of the
Riemann hypothesis, the standard conjectures retain their interest for the
theory of motives.

The first, the {\it Lefschetz standard conjecture} (Grothendieck 1969, \S
3), states that, for a smooth projective variety $V$ over an algebraically
closed field, the operators $\Lambda$, rendering commutative the diagrams ($%
0\leq r\leq 2n$, $n=\dim V)$ 
$$
\begin{array}{ccc} 
H^{r}(V)&  {{L^{n-r}}\atop{\lrar\atop\approx}}&H^{2n-r}(V)\\ 
\big\downarrow{\scriptstyle\Lambda}&&\big\downarrow{\scriptstyle L} \\[4pt]
H^{r-2}(V) & {L^{n-r+2}\atop{\llrar\atop\approx}}&H^{2n-r+2}(V)\; ,\end{array}
$$  
are algebraic. Here $H$ is a Weil cohomology theory and $L$ is cup product
with the class of a smooth hyperplane section ($L^{n-r}$ is assumed to be an
isomorphism for $n\geq r$, and $L^{n-r}=(L^{r-n})^{-1}$ for $n<r$). This
conjecture is known for abelian varieties (Lieberman 1968, Kleiman 1968),
surfaces and Weil cohomologies for which $\dim H^{1}(V)=2\dim \Pic^{0}(V)$
(Grothendieck), and a few other varieties (see Kleiman 1994, 4.3). For
abelian varieties, it is even known that the operator $\Lambda $ is defined
by a Lefschetz class, i.e., a class in the ${\Bbb Q}{}$-algebra generated
by divisor classes (Milne 1999a, 5.9).

The second, the {\it Hodge standard conjecture} (Grothendieck 1969, \S 4),
states that, for $r\leq n/2$, the bilinear form%
$$
(x,y)\mapsto (-1)^{r}\langle L^{n-2r}x\cdot y\rangle \colon P^{r}(V)\times
P^{r}(V)\rightarrow {\Bbb Q}
$$%
is positive-definite. Here $P^{r}(V)$ is the ${\Bbb Q}{}$-space of
primitive algebraic classes of codimension $r$ modulo homological
equivalence. In characteristic zero, Hdg($V$) is a consequence of Hodge
theory (Weil 1958). In nonzero characteristic, Hdg($V$) is known for
surfaces (Segre 1937; Grothendieck 1958). An important consequence of the
Hodge standard conjecture for abelian varieties, namely, the positivity of
the Rosati involution, was proved in nonzero characteristic by Weil (1948, Th\'eor\`eme 38). Apart from these
examples and the general coherence of Grothendieck's vision, there appears to have been little evidence for the
conjecture in nonzero characteristic.

In this paper, we prove that the Hodge standard conjecture holds for abelian
varieties in arbitrary characteristic if the Hodge conjecture holds for
complex abelian varieties of CM-type. 

Let $\Mot({\Bbb F};{\cal  A}{})$ be the category of motives based on
abelian varieties over ${\Bbb F}{}$ using the numerical equivalence
classes of algebraic cycles as correspondences. This is a Tannakian category
(Jannsen 1992, Deligne 1990), and it is known that the Tate conjecture for
abelian varieties over finite fields implies that it has all the major
expected properties but one, namely, that the Weil forms coming from
algebraic geometry are positive for the canonical polarization on $\Mot(%
{\Bbb F};{\cal  A}{})$ (see Milne 1994, especially 2.47).

In Milne 1999b it is shown that the Hodge conjecture for complex abelian
varieties of CM-type is stronger than (that is, implies) the Tate conjecture
for abelian varieties over finite fields. Here, we show that the stronger
conjecture also implies the positivity of the Weil forms coming from
algebraic geometry (Theorem \ref{f}). As a consequence, we obtain the Hodge
standard conjecture for abelian varieties over finite fields, and a
specialization argument then proves it over any field of nonzero
characteristic (Theorem \ref{r}).

Most of the arguments in the paper hold with ``algebraic cycle'' replaced by
``Lefschetz cycle''. In fact, the analogue of the Hodge standard conjecture
holds unconditionally for Lefschetz classes on abelian varieties. In
particular, the Hodge standard conjecture is true for abelian varieties
without exotic (i.e., non-Lefschetz) algebraic classes\ (\ref{t5}, \ref{t7}).

In preparation for proving these results, we study in Section 1 the polarizations
on a quotient Tannakian category.

\demo{Notation  and Conventions}
The algebraic closure of ${\Bbb Q}{}$ in ${\Bbb C}{}$ is denoted~${\Bbb Q}{}^{{\rm al} }$. We fix a $p$-adic prime on
${\Bbb Q}{}^{{\rm  al}}$ and denote its residue field by ${\Bbb F}{}$.

By the Hodge conjecture for a variety $V$ over ${\Bbb C}{}$, we mean the
statement that, for all $r$, the ${\Bbb Q}{}$-space $H^{2r}(V,{\Bbb Q}%
{})\cap H^{r,r}$ is spanned by the classes of algebraic cycles.

By the Tate conjecture for a variety $V$ over a finite field ${\Bbb F}%
{}_{q}$ we mean the statement that, for all $r$, the order of the pole of
the zeta function $Z(V,t)$ at $t=q^{-r}$ is equal to the rank of the group
of numerical equivalence classes of algebraic cycles of codimension $r$ on $%
V $ (Tate 1994, 2.9). We say that a variety over ${\Bbb F}{}$ satisfies
the Tate conjecture if all of its models over finite fields satisfy the Tate
conjecture (equivalently, one model over a ``sufficiently large'' finite
field).

For abelian varieties $A$ and $B$, $\Hom(A,B)_{{\Bbb Q}{}}=\Hom%
(A,B)\otimes {\Bbb Q}{}$. An abelian variety $A$ over ${\Bbb C}{}$ (or $%
{\Bbb Q}{}^{{\rm al} }$) is said to be of CM-type if, for each simple
isogeny factor $B$ of $A$, $\End(B)_{{\Bbb Q}{}}$ is a commutative field
of degree $2\dim B$ over ${\Bbb Q}{}$. A polarization of $A$ is the
isogeny $A\rightarrow A^{\vee }$ from $A$ to its dual defined by an ample
divisor on $A$.

Let ${\cal  S}{}$ be a set of smooth projective varieties over an
algebraically closed field $k$ satisfying the following condition:

\begin{description}
\item[(0.1)] the projective spaces are in ${\cal  S}{}$ and ${\cal  S}{}$
is closed under passage to a connected component and under the formation of
products and disjoint unions.
\end{description}

\noindent For example, ${\cal  S}{}$ could be the class ${\cal  T}{}$ of
all smooth projective varieties over $k$ or the smallest class ${\cal  A}%
{} $ satisfying (0.1) and containing the abelian varieties. \noindent Then $%
\Mot(k;{\cal  S})$ is defined to be the category of motives based on the
abelian varieties over $k$ with the algebraic classes modulo numerical
equivalence as correspondences.
\enddemo

{\it Acknowledgement}.
I thank the referee for his suggestions for simplifying and shortening the
article.

\section{Polarizations on quotient categories} 

We refer the reader to Deligne 1989, \S \S 5, 6,   for the theory of algebraic
geometry in a Tannakian category ${\bf C}{}$. In particular, the
fundamental group{\it \ }$\pi ({\bf C})$ of ${\bf C}{}$ is an affine
group scheme in ${\bf C}{}$, such that, for any fibre functor $\omega $
on ${\bf C}{}$,%
$$
\underline{\Aut}^{\otimes }(\omega )\cong \omega (\pi {\bf (C}{})).
$$%
The fundamental group acts on the objects of ${\bf C}{}$. When $H$ is a
closed subgroup of $\pi ({\bf C}{})$, we let $X^{H}$ denote the largest
subobject of $X$ on which the action of $H$ is trivial, and we let ${\bf C%
}{}^{H}$ denote the full subcategory of ${\bf C}{}$ of objects on which
the action of $H$ is trivial. The functor $\Hom(\1,-)$ is a tensor
equivalence from ${\bf C}{}^{\pi ({\bf C}{})}$ to the category of
finite-dimensional vector spaces over $F=_{{\rm df}}\End(\1)$, which allows
us to regard the objects of ${\bf C}{}^{\pi ({\bf C}{})}$ as vector
spaces. When $\pi ({\bf C}{})$ is commutative, it lies in $\Ind({\bf C}%
{}^{\pi ({\bf C}{})})$, and hence can be regarded as a group scheme in
the usual sense.

We refer to Saavedra 1972, V 2.3.1, V 3.2.1,  for the definitions of a Weil
form and of a (graded) polarization on a Tate triple over ${\Bbb R}{}$. We
define a polarization on a Tate triple ${\bf C}{}$ over ${\Bbb Q}{}$ to
be a polarization on ${\bf C}_{({\Bbb R}{})}$.

\numbereddemo{Remark} 
\label{a}
Let $({\bf C}{},w,{\Bbb T)}$ be a Tate triple. In particular, $%
{\bf C}{}$ is a rigid tensor category, and so each object $X$ has a dual $%
(X^{\vee },\ev_{X})$; moreover,%
\begin{equation}
\underline{\End}(X)^{\vee }\cong (X^{\vee }\otimes X)^{\vee }\cong X\otimes
X^{\vee }\cong \underline{\End}(X^{\vee }).\hskip.75in \speqnu{1.1.1} \label{a1}
\end{equation}%
Let $X$ be an object of weight $n$ in ${\bf C}{}$. A nondegenerate $%
(-1)^{n}$-symmetric bilinear form $\psi \colon X\otimes X\rightarrow \1(-n)$
on $X$ defines an isomorphism $X\rightarrow X^{\vee }(-n)$, and hence an
isomorphism%
$$
\underline{\End}(X)\rightarrow \underline{\End}(X^{\vee }(-n))\cong 
\underline{\End}(X^{\vee }).
$$%
This, together with the pairing%
$$
\ev\colon \underline{\End}(X)^{\vee }\otimes \underline{\End}(X)\rightarrow %
\1
$$%
and the isomorphism (\ref{a1}), gives a symmetric bilinear form%
$$
T^{\psi }\colon \underline{\End}(X)\otimes \underline{\End}(X)\rightarrow \1.
$$%
On $\End(X)\subset \underline{\End}(X)$, $T^{\psi }$ is the form $%
(u,v)\mapsto \Tr_{X}(u^{\psi }\cdot v)$, and so to say that $\psi $ is a
Weil form amounts to saying that the form induced by $T^{\psi }$ on $\End(X)$
is positive-definite.
\enddemo

A morphism $F\colon ({\bf C}{}_{1},w_{1},{\Bbb T}_{1})\rightarrow (%
{\bf C}{}_{2},w_{2},{\Bbb T}_{2})$ of Tate triples is an exact tensor
functor $F\colon {\bf C}{}_{1}\rightarrow {\bf C}{}_{2}$ preserving
the gradations together with an isomorphism $F({\Bbb T}_{1})\cong {\Bbb %
{\Bbb T}}_{2}$. We say that such a morphism $F$ maps a polarization ${\mathnormal{\Pi}}
_{1}$ on ${\bf C}{}_{1}$ to a polarization ${\mathnormal{\Pi}}_{2}$ on ${\bf C}%
{}_{2} $ (denoted $F\colon {\mathnormal{\Pi}}_{1}\mapsto {\mathnormal{\Pi}}_{2}$) if 
$$
\psi \in {\mathnormal{\Pi}}_{1}(X)\Rightarrow F\psi \in {\mathnormal{\Pi}}_{2}(FX),
$$%
in which case, for an $X$ of weight $n$, ${\mathnormal{\Pi}}_{1}(X)$ consists of the
bilinear forms $\psi \colon X\otimes X\rightarrow \1(-n)$ such that $F\psi
\in {\mathnormal{\Pi}}_{2}(FX)$. In particular, given $F$ and ${\mathnormal{\Pi}}_{2}$, there exists at
most one polarization ${\mathnormal{\Pi}}_{1}$ on ${\bf T}_{1}$ such that $F\colon {\mathnormal{\Pi}}_{1}\mapsto $ ${\mathnormal{\Pi}}_{2}$.

\proclaim{Lemma}
\label{b}
Let $F\colon {\bf C}{}\rightarrow {\bf Q}{}$ be a morphism of
Tate triples{\rm ,} and assume that every object of ${\bf Q}{}$ is a direct
summand of an object in the image of $F${\rm .} Let ${\mathnormal{\Pi}}_{1}$ be a polarization
on ${\bf C}{}${\rm .} There exists a polarization ${\mathnormal{\Pi}}_{2}$ on ${\bf Q}{}$
such that $F\colon {\mathnormal{\Pi}}_{1}\mapsto {\mathnormal{\Pi}}_{2}$ if and only if{\rm ,} for all $X$ in $%
{\bf C}{}$ and all $\psi \in {\mathnormal{\Pi}}_{1}(X)${\rm ,} $F\psi $ is a Weil form on $FX${\rm .}
\endproclaim

\demo{Proof}
$\Rightarrow $: This follows directly from the definitions.

$\Leftarrow $: For $Y$ a direct summand of $FX$, define ${\mathnormal{\Pi}}_{2}(Y)$ to be
the compatibility class of $(F\psi )|Y$ for some $\psi \in {\mathnormal{\Pi}}_{1}(X)$. It
is straightforward to verify that the sets ${\mathnormal{\Pi}}_{2}(Y)$ are well-defined
and form a polarization on ${\bf Q}{}$.
\enddemo

Recall that an exact tensor functor $q\colon {\bf C}{}\rightarrow {\bf %
Q}{}$ of Tannakian categories defines a morphism $\pi (q)\colon \pi ({\bf %
Q}{})\rightarrow q(\pi ({\bf C}{}))$ (Deligne 1990, 8.15.2).

\numbereddemo{Definition}
\label{c}Let $q\colon {\bf C}{}\rightarrow {\bf Q}{}$ be an exact
tensor functor, and let $H$ be a closed subgroup of ${\bf \pi }({\bf C}%
{})$. We say that $({\bf Q}{},q)$ is a {\it quotient of }${\bf C}{}$ 
{\it by }$H$ if $\pi (q)$ is an isomorphism of $\pi {\bf (Q}{})$ onto $%
q(H)$.
\enddemo

When $({\bf Q}{},q)$ is a quotient of ${\bf C}{}$ by $H\subset \pi (%
{\bf C}{})$, every object in ${\bf Q}{}$ is a subquotient of an object
in the image of $q$. Moreover, $q$ maps ${\bf C}{}^{H}$ into ${\bf Q}%
{}^{\pi ({\bf Q}{})}$, and so, for $X\in {\bf C}{}^{H}$, we can
identify $qX$ with the vector space $\Hom(\1,qX)$. With this identification,
there is a functorial isomorphism 
$$
\Hom_{{\bf Q}{}}(qX,qY)\cong q(\underline{\Hom}(X,Y)^{H}),X,Y\in %
\ob({\bf C}{}).
$$

\proclaim{Proposition}
\label{d}
Let $({\bf C}{},w,{\Bbb T)}$ be a Tate triple over ${\Bbb R}%
{}${\rm .} Let $({\bf Q,}q)$ be a quotient of ${\bf C}{}{}$ by $H\subset \pi
({\bf C}{})${\rm ,} and let ${\mathnormal{\Pi}}_{1}$ be a polarization on ${\bf C}{}${\rm .}
Suppose that $H\supset w({\Bbb G}_{m})${\rm ,} so that ${\bf Q}{}$ inherits a
Tate triple structure from that on ${\bf C}{}${\rm ,} and that ${\bf Q}{}$
is semisimple{\rm .} Assume\/{\rm :}\/

\begin{itemize}
\item[{\rm (*)}] for all $X$ in ${\bf C}{}^{H}$ and all $\psi \in {\mathnormal{\Pi}}_{1}(X)${\rm ,} 
$q\psi $ is a positive\/{\rm -}\/definite form on the vector space $qX${\rm .}
\end{itemize}

\noindent Then there exists a polarizaton ${\mathnormal{\Pi}}_{2}$ on ${\bf Q}{}$
such that $q\colon {\mathnormal{\Pi}}_{1}\mapsto {\mathnormal{\Pi}}_{2}${\rm .}
\endproclaim 

\demo{Proof}
Because ${\bf Q}{}$ is semisimple, every object of ${\bf Q}{}$ is a
direct summand of an object in the image of $q$. We shall check the
condition in Lemma \ref{b}.

Let $\psi \in {\mathnormal{\Pi}}_{1}(X)$. Then $T^{\psi }$ is positive for ${\mathnormal{\Pi}}_{1}$, and
hence so also is its restriction $T^{\psi }|$ to $\underline{\End}(X)^{H}$.
Therefore, (*) implies that $q(T^{\psi }|)$ is a positive-definite form on
the vector space $q(\underline{\End}(X)^{H})$. But $q(\underline{\End}%
(X)^{H})\cong \End_{{\bf Q}{}}(qX)$ and $q(T^{\psi }|)\cong T^{q\psi }$,
and so $q\psi $ is a Weil form, as required.
\enddemo

\numbereddemo{Remark} 
\label{e}
Instead of (*), it suffices to assume that there exists a single $X$
in ${\bf C}{}^{H}$ such that $\pi ({\bf C}{})/H$ acts faithfully on $X$
and a single $\psi \in \Pi (X)$ such that $q\psi $ is a positive-definite
form on $qX$.
\enddemo

\vglue-12pt
\section{Polarizations on categories of motives over finite fields}
\vglue-8pt
Consider $\Mot(k;{\cal  S}{})$ for ${\cal  S}{}$ some class satisfying
(0.1). For an abelian variety $A$ in ${\cal  S}{}$, a divisor $D$ on $A$
defines a pairing $\phi _{D}\colon h_{1}A\otimes h_{1}A\rightarrow {\Bbb %
{\Bbb T}}$, which is a Weil form if $D$ is ample (Weil 1948, Th\'eor\`eme 38). Such a Weil form will be said to be {\it
geometric}.

Consider $\Mot({\Bbb F}{};{\cal  A}{})$. If the Tate conjecture holds
for all abelian varieties over ${\Bbb F}{}$, then $\Mot({\Bbb F}{};%
{\cal  A}{})$ is a semisimple Tate triple over ${\Bbb Q}$ with the Weil
number torus $P$ as fundamental group (see, for example, Milne 1994, 2.26).
Moreover, there exist two graded polarizations on $\Mot({\Bbb F})$, and
for exactly one of these (denoted $\Pi ^{{\rm Mot} }$) the geometric Weil
forms on any supersingular elliptic curve are positive (ibid., 2.44).

Consider $\Mot({\Bbb Q}^{{\rm al} };{\cal  C}{})$ where ${\cal  C}{}$
is the smallest class satisfying (0.1) and containing the abelian varieties
of CM-type over ${\Bbb Q}{}^{{\rm al} }$. It is a Tate triple over $%
{\Bbb Q}{}$ with the Serre group $S$ as fundamental group, and it has a
canonical polarization $\Pi ^{{\rm CM} }$. If the Hodge conjecture holds for
complex abelian varieties of CM-type, then the Tate conjecture holds for
abelian varieties over ${\Bbb F}{}$ (Milne 1999b, 7.1), and, corresponding
to the $p$-adic prime we have fixed on ${\Bbb Q}{}^{{\rm al} }$, there is
a reduction functor $R\colon \Mot({\Bbb Q}{}^{{\rm al} };{\cal  C}%
{})\rightarrow \Mot({\Bbb F};{\cal  A}{})$ which realizes $\Mot({\Bbb F%
}{}{};{\cal  A}{})$ as the quotient of $\Mot({\Bbb Q}{}^{{\rm al} },%
{\cal  C})$ by the closed subgroup $P$ of the Serre group $S$ (a
description of the inclusion $P\hookrightarrow S$ can be found, for example,
in Milne 1994, 4.12). For a motive $X=h(A)(r)$ in $\Mot({\Bbb Q}^{{\rm al} %
};{\cal  C}{})$, $R(X^{P})$ is the ${\Bbb Q}{}$-space of numerical
equivalence classes of algebraic cycles of codimension $r$ on the reduction $%
A_{{\Bbb F}{}}$ of~$A$.

\proclaim{Theorem}
\label{f}
If the Hodge conjecture holds for complex abelian varieties of
{\rm CM-}\/type{\rm ,} then $R\colon \Pi ^{{\rm CM} }\mapsto \Pi ^{{\rm Mot} }$ and all
geometric Weil forms on all abelian varieties over ${\Bbb F}{}$ are
positive for $\Pi ^{{\rm Mot} }${\rm .}
\endproclaim

\demo{Proof}
I claim that to prove the theorem it suffices to show:

 \begin{quote}%
 (*) there exists a polarization \negthinspace $\Pi $ \negthinspace on
\negthinspace $\Mot({\Bbb F})$ such that $R\colon \Pi ^{{\rm CM} %
}\!\mapsto \!\Pi $.
 \end{quote} 

\noindent The geometric Weil forms are positive for $\Pi ^{{\rm CM} }$%
\noindent\ and every polarized abelian variety $A$ over ${\Bbb F}{}$ lifts
(up to isogeny) to a polarized abelian variety of CM-type over ${\Bbb Q}%
{}^{{\rm al} }$ (Zink 1983, 2.7) and so if $R\colon \Pi ^{{\rm CM} }\mapsto
\Pi $, then every geometric Weil form is positive for $\Pi $. In particular,
the geometric Weil forms on a supersingular elliptic curve are positive, and
so $\Pi =\Pi ^{{\rm Mot} }$. This proves the claim.

We now prove (*). Fix a CM-field $K\subset {\Bbb {\Bbb Q}{}}^{{\rm al} }$
such that

\begin{itemize}
\item[(a)] $K$ is finite and Galois over ${\Bbb Q}{}$, and

\item[(b)] $K$ properly contains an imaginary quadratic field in which $p$ splits.
\end{itemize}

\noindent Let ${\cal  C}{}^{K}$ (resp.\ ${\cal  A}{}^{K}$) be the
smallest subset of ${\cal  C}{}$ (resp.\ ${\cal  A}{}$) satisfying (0.1)
and containing the CM abelian varieties over ${\Bbb Q}{}^{{\rm al} }$ with
reflex field contained in $K$ (resp.\ the abelian varieties over ${\Bbb F}%
{} $ with endomorphism algebra split by $K$), and let $S^{K}$ and $P^{K}$ be
the corresponding quotients of $S$ and $P$. It suffices to prove (*) for 
$$
R^{K}\colon \Mot({\Bbb Q}^{{\rm al} };{\cal  {\cal  C}{}}%
^{K})\rightarrow \Mot({\Bbb F}{};{\cal  A}^{K}).
$$

Let $A$ be the product of a set of representatives for the simple isogeny
classes of abelian varieties in ${\cal  C}{}^{K}$, and let $X=\underline{%
\End}(h_{1}A)^{P}$. It follows from Milne 1999b that $S^{K}/P^{K}$ acts
faithfully on $X$: with the notation  of that paper, $T^{K}$ acts faithfully
on $h_{1}A$ as an object of the category of Lefschetz motives; therefore, $%
T^{K}/L^{K}$ acts faithfully on $\underline{\End}(h_{1}A)^{L^{K}}$, and
(ibid.\ \S 6) the canonical map 
$$
S^{K}/P^{K}\rightarrow T^{K}/L^{K}
$$%
is injective.

Let $\phi $ be the geometric Weil form on $h_{1}A$ defined by an ample
divisor $D$, and let $\psi =$ $T^{\phi }|X$. Then $\psi \in $ $\Pi ^{{\rm CM}}(X)$, and it suffices to show that
$R^{K}(\psi )$ is positive-definite (%
\ref{d}, \ref{e}). But $R^{K}(X)=\End(A_{{\Bbb F}{}})_{{\Bbb Q}{}}$ and $%
R^{K}(\phi )$ is the trace pairing $u,v\mapsto \Tr(u\cdot v^{\dagger })$ of
the Rosati involution defined by $D_{{\Bbb F}{}}$, which is
positive-definite (Weil 1948, Th\'eor\`eme 38).
\enddemo

\section{The Hodge standard conjecture}

Throughout this section, $k$ is an algebraically closed field and ${\cal  S%
}{}$ is a class of smooth projective varieties over $k{}$ satisfying (0.1).

By a Weil cohomology theory on ${\cal  S}{}$, we mean a contravariant
functor $V\mapsto H^{\ast }(V)$ satisfying the conditions (1)--(4), (6) of
Kleiman  1994, \S 3  (finiteness, Poincar\'{e} duality, K\"{u}nneth formula,
cycle map, strong Lefschetz theorem), except that we remember the Tate
twists. For example, $\ell $-adic \'{e}tale cohomology, $\ell \neq \ch(k)$,
is a Weil cohomology theory in this sense (the strong Lefschetz theorem is
proved in Deligne 1980).

For $V\in {\cal  S}{}$, $A_{\sim }^{\ast }(V)$ denotes the ${\Bbb Q}{}$%
-algebra of algebraic classes on $V$ modulo an admissible equivalence
relation $\sim $, for example, numerical equivalence (num), or homological
equivalence (hom) with respect to a Weil cohomology~$H$.

We say that a Weil cohomology theory is {\it good} if homological
equivalence coincides with numerical equivalence on algebraic cycles with $%
{\Bbb Q}$-coefficients for all varieties in ${\cal  S}{}$.

Let $H$ be a Weil cohomology theory on ${\cal  S}{}$. For a connected
variety $V$ in ${\cal  S}{}$ of dimension $n$, define $P^{r}(V)$ to be the
subspace of $A_{{\rm hom} }^{r}(V)$ on which $L^{n-2r+1}$ is zero. Let $%
\theta ^{r}$ be the bilinear form%
$$
(x,y)\mapsto (-1)^{r}\langle L^{n-2r}x\cdot y\rangle \colon P^{r}(V)\times
P^{r}(V)\rightarrow {\Bbb Q}{},\quad r\leq n/2.
$$%
As originally stated  (Grothendieck 1969), the Hodge standard conjecture
asserts that these pairings are positive-definite when $H$ is $\ell $-adic 
\'{e}tale cohomology. Kleiman (1968, \S3) states the conjecture for any
Weil cohomology theory. Note that if the Hodge standard conjecture holds for
one good Weil cohomology theory, then it holds for all.

For a Weil cohomology theory $H$, $\pi ^{r}$ denotes the projection onto $%
H^{r}$ and $\Lambda $, $^{c}\Lambda $, $\ast $, $p^{r}$ denote the maps
defined in Kleiman 1968, 1.4 (corrected in Kleiman 1994, \S 4).

\proclaim{Proposition}
\label{n}For all good Weil cohomology theories $H$ on ${\cal  S}{}${\rm ,} the
operators $\Lambda ${\rm ,} $^{c}\Lambda ${\rm ,} $\ast ${\rm ,} $p^{r}${\rm ,} and $\pi ^{r}$ are
defined by algebraic cycles that{\rm ,} modulo numerical equivalence{\rm ,} depend only
on $L$ {\rm (}\/not $H${\rm ).}
\endproclaim 

\demo{Proof}
Let $H$ be a good Weil cohomology theory on ${\cal  S}{}$. By assumption, the
Lefschetz standard conjecture holds for all $V\in {\cal  S}{}$,  and the proposition can be proved as in Kleiman 1994,
5.4  (the Hodge standard conjecture is used there only to deduce that numerical
equivalence coincides with homological equivalence on $V\times V$).
\enddemo

When there exists a good Weil cohomology theory on ${\cal  S}{}$, we
define\break $\Mot(k{};{\cal  S}{})$ to be the category of motives with the
commutativity constraint modified using the $\pi ^{r}$'s given by (\ref{n}).
It is semisimple (Jannsen 1992), hence Tannakian (Deligne 1990), and it has
a natural structure of a Tate triple. Let $V\in {\cal  S}{}$ be connected 
of dimension $n$, and let $Z$ be a smooth hyperplane section of $V$. Then $%
l=_{{\rm df}}\Delta _{V}(Z)\in A_{{\rm num} }^{n+1}(V\times V)$ is a
morphism 
$$
l\colon h(V)\rightarrow h(V)(1)
$$%
of degree $2$. Define%
$$
\varphi ^{r}\colon h^{r}(V)\otimes h^{r}(V)\rightarrow \1(-r)
$$%
to be the composite%
$$
h^{r}(V)\otimes h^{r}(V)\stackrel{{\rm id}\otimes \ast }{\longrightarrow }%
h^{r}(V)\otimes h^{2n-r}(V)(n-r)\rightarrow h^{2n}(V)(n-r)\cong \1(-r)
$$%
($\ast $ as in \ref{n}). Let $p{}^{r}(V)$ be the largest subobject of 
$$
\Ker(l^{n-2r+1}\colon h^{2r}(V)(r)\rightarrow h^{2n-2r+2}(V)(n-r+1))
$$%
on which $\pi =_{{\rm df}}\pi (\Mot(k;{\cal  S}{}))$ acts trivially. For
any good Weil cohomology theory $H$ on ${\cal  S}{}$, 
$$
\gamma (p^{r}(V))\cong P^{r}(V)
$$%
where $\gamma $ is the tensor equivalence $\Hom(\1,-)$ from $\Mot(k;{\cal  %
S}{})^{\pi }$ to finite-dimensional ${\Bbb Q}{}$-vector spaces, and there
is a pairing 
$$
\vartheta ^{r}\colon p^{r}(V)\otimes p^{r}(V)\rightarrow \1,
$$%
also fixed by $\pi $, such that $\gamma (\vartheta ^{r})=$ $\theta ^{r}$.

\proclaim{Proposition}
\label{q}
Let $H$ be a good Weil cohomology theory on ${\cal  S}{}${\rm .} The
following statements are equivalent\/{\rm :}

\begin{itemize}
\item[{\rm (a)}] the Hodge standard conjecture holds for $H$ and the varieties in $%
{\cal  S}{}${\rm ;}

\item[{\rm (b)}] there exists a polarization $\Pi $ on $\Mot(k;{\cal  S}{})$ for
which the forms $\varphi ^{r}$ are positive\/{\rm ;}\/

\item[{\rm (c)}] there exists a polarization $\Pi $ on $\Mot(k;{\cal  S}{})$ for
which the forms $\vartheta ^{r}$ are positive{\rm .}
\end{itemize}

\endproclaim 

\vglue-4pt
{\it Proof}.
(a)$\Rightarrow $(b): See Saavedra 1972, VI 4.4.

(b)$\Rightarrow $(c): The restriction of $\varphi ^{2r}\otimes {\rm id}_{\1(2r)}$
to the subobject $p^{r}(V)$ of $h^{2r}(V)(r)$ is the form $\vartheta ^{r}$.

(c)$\Rightarrow $(a): Let $\Pi $ be a polarization on $\Mot(k;{\cal  S}{})$
for which the forms $\vartheta ^{r}$ are positive. The restriction of $\Pi $
to $\Mot(k{};{\cal  S}{})^{w({\Bbb G}_{m})}$ is a symmetric
polarization, and so there exists an ${\Bbb R}{}$-valued fibre functor $%
\omega $ on $\Mot(k{};{\cal  S}{})^{w({\Bbb G}_{m})}$ carrying $\Pi $%
-positive forms to positive-definite symmetric forms (Deligne and Milne
1982, 4.27). The restriction of $\omega $ to $\Mot(k;{\cal  S}{})_{(%
{\Bbb R}{})}^{\pi }$ is (uniquely) isomorphic to $\gamma $, and so $\gamma
(\vartheta ^{r})$ is positive-definite.
\hfill\qed

\proclaim{Theorem}
\label{r}
Let $k$ be an algebraically closed field{\rm .} If the Hodge conjecture
holds for complex abelian varieties of  {\rm CM-}\/type{\rm ,} then{\rm ,} for all $\ell \neq \ch%
(k)${\rm ,}
\begin{itemize}
\item[{\rm (a)}] numerical equivalence coincides with $\ell $\/{\rm -}\/adic {\rm \'{\it e}}tale
homological equivalence on abelian varieties over $k${\rm ,} and

\item[{\rm (b)}] the Hodge standard conjecture holds for all abelian varieties over $k$
and the $\ell $\/{\rm -}\/adic {\rm \'{\it e}}tale cohomology theory{\rm .}
\end{itemize}

\endproclaim

\demo{Proof}
 (a)  {\it for}  $k={\Bbb F}{}$.{\it \ }The Hodge conjecture for complex
abelian varieties of CM-type implies the Tate conjecture (Milne 1999b, 7.1),
which implies (a) (see, for example, Tate 1994, 2.7).

 (b)   {\it for}  $k={\Bbb F}{}$. For abelian varieties over ${\Bbb Q}{}^{%
{\rm al} }$, the Betti cohomology theory is good (Lieberman 1968) and the
Hodge standard conjecture holds, and so by (\ref{q}) there is a polarization $%
\Pi $ on $\Mot({\Bbb Q}^{{\rm al} };{\cal  {\cal  C}{}})$ for which
the forms 
$$
\varphi ^{r}\colon h^{r}(A)\otimes h^{r}(A)\rightarrow \1(-r)
$$%
are positive. Clearly, $\Pi $ is the canonical polarization $\Pi ^{{\rm CM} }
$ in Section 2. Let $Z$ be the hyperplane section of $A$ used in the definition
of $\varphi ^{r}$. Because $R\colon \Pi ^{{\rm CM} }\mapsto \Pi ^{{\rm Mot} }
$ (Theorem \ref{f}), the form%
$$
\varphi ^{r}\colon h^{r}(A_{{\Bbb F}{}})\otimes h^{r}(A_{{\Bbb F}%
{}})\rightarrow \1(-r)
$$%
defined by the reduction $Z_{{\Bbb F}{}}$ of $Z$ on $A_{{\Bbb F}{}}$ is
positive for $\Pi ^{{\rm Mot} }$. Every polarized abelian variety $A$ over $%
{\Bbb F}{}$ lifts (up to isogeny) to a polarized abelian variety of
CM-type over ${\Bbb Q}{}^{{\rm al} }$ (Zink 1983, 2.7), and so 
$$
(A_{{\Bbb F}{}},Z_{{\Bbb F}{}}\hbox{ modulo numerical equivalence})
$$%
is arbitrary. Proposition \ref{q} now gives (b).

 (a, b) {\it for arbitrary}  $k$. For an abelian variety $A$ of dimension $n$
over $k$, consider the commutative diagram:%
\figin{fig1}{1000}%
\noindent There is a similar diagram for a smooth specialization $A_{%
{\Bbb F}{}}$ of $A$ to an abelian variety over ${\Bbb F}{}$. The
specialization maps on the cohomology groups are bijective and hence they
are injective on the $P$'s. Since the pairings are compatible, this implies
the Hodge standard conjecture for $A$ and $\ell $-adic \'{e}tale cohomology.
Because the Lefschetz standard conjecture is known for abelian varieties,
this in turn implies that numerical equivalence coincides with $\ell $-adic
homological equivalence for $A$ (Kleiman 1994, 5-4).
\enddemo 

\proclaim{{C}orollary}
\label{t} If the Hodge conjecture holds for complex abelian varieties of
{\rm CM-}\/type{\rm ,} then{\rm ,} for any algebraically closed field $k${\rm ,} $\Mot(k;{\cal  A}%
{}) $ has a polarization {\rm (}\/necessarily unique\/{\rm )} for which the forms $\vartheta
^{r} $ and $\varphi ^{r}$ are positive{\rm .}
\endproclaim

\demo{Proof}
For $\ell \neq \ch(k)$, Theorem \ref{r}(a) shows that $\ell $-adic \'{e}tale
cohomology is good. Now apply (\ref{r}b) and (\ref{q}).
\enddemo

\numbereddemo{Remark} 
\label{u}Assume the Hodge conjecture holds for complex abelian varieties of
CM-type, and let $H$ be a Weil cohomology theory on ${\cal  A}{}$ (over an
algebraically closed field $k$). Because the Lefschetz standard conjecture
is known for abelian varieties, if $H$ is not good, then the Hodge standard
conjecture fails for $H$ (Kleiman 1994, 5-1). Thus, the Hodge standard
conjecture holds for $H$ if and only if $H$ is good.
\enddemo

\numbereddemo{Remark} 
\label{s4}Let $K$ be a CM-subfield of ${\Bbb Q}{}^{{\rm al} }$ satisfying
conditions (a) and (b) of the proof of (\ref{f}). The preceding arguments
can be modified to show that, if the Hodge conjecture holds for all complex
abelian varieties with reflex field contained in $K$, then the conclusions
of Theorem \ref{r} hold for each abelian variety over ${\Bbb F}{}$ whose
endomorphism algebra is split by $K$. In fact, condition (b) is not
necessary for this statement because, as Deligne pointed out to me, the
results of Milne 1999b, \S 6   hold without it.
\enddemo

\numbereddemo{Remark} 
\label{t5}Most of the preceding arguments hold with ``algebraic cycle''
replaced by ``Lefschetz cycle'' (cf.\ Milne 1999a, \S 5). Let $A$ be an
abelian variety over $k$. Recall that, for any Weil cohomology theory, if a
Lefschetz class $a$ on $A$ is not homologically equivalent to zero, then
there exists a Lefschetz class $b$ on $A$ of complementary dimension such
that $\langle a\cdot b\rangle \neq 0$ (ibid.\ 5.2). Thus, homological
equivalence on Lefschetz classses is independent of the Weil cohomology
theory, and coincides with numerical equivalence.

Let $D^{r}(A)$ be the ${\Bbb Q}{}$-space of Lefschetz classes modulo
numerical equivalence on $A$ of codimension $r$, and let $DP^{r}(A)$ be the $%
{\Bbb Q}{}$-subspace on which $L^{n-2r+1}$ is zero. With the notations of
Milne 1999b, the categories $\LCM({\Bbb Q}{}^{{\rm al} })$ and $\LMot(%
{\Bbb F})$ of Lefschetz motives have canonical polarizations, and the
reduction functor $\LCM({\Bbb Q}{}^{{\rm al} })\rightarrow \LMot({\Bbb F}%
)$ maps one to the other. The same argument as in the proof of Theorem \ref%
{r} shows that the bilinear forms%
$$
(x,y)\mapsto (-1)^{r}\langle L^{n-2r}x\cdot y\rangle \colon DP^{r}(A)\times
DP^{r}(A)\rightarrow {\Bbb Q}
$$%
are positive-definite for $r\leq n/2$ and $A$ an abelian variety over $%
{\Bbb F}{}$. In other words, the Lefschetz analogue of the Hodge standard
conjecture holds unconditionally for abelian varieties over ${\Bbb F}{}$.
A specialization argument (as in the proof of Theorem~\ref{r}) extends the statement
to arbitrary $k$.
\enddemo

\numbereddemo{Remark} 
\label{t7}Recall that a Hodge, Tate, or algebraic class on a variety is said
to be {\it exotic }if it is not Lefschetz. Remark \ref{t5} shows that the
Hodge standard conjecture holds unconditionally for abelian varieties with
no exotic algebraic classes. For examples (discovered by Lenstra, Spiess,
and Zarhin) of abelian varieties over ${\Bbb F}{}$ with no exotic Tate
classes, and hence no exotic algebraic classes, see Milne 2001, A.7.
\enddemo

\numbereddemo{Remark} 
\label{v}Grothendieck (1969) stated: ``Alongside the problem of resolution
of singularities, the proof of the standard conjectures seems to me to be
the most urgent task in algebraic geometry.'' Should the Hodge conjecture
remain inaccessible, even for abelian varieties of CM-type, Theorem \ref{r}\
suggests a possible approach to proving the Hodge standard conjecture for
abelian varieties, namely, improve the theory of absolute Hodge classes
(Deligne 1982) sufficiently to remove the hypothesis from the theorem.
\enddemo

% {\centerline{\eightsc References} }
  
\def\joinrel{\mathrel{\mkern-8mu}}
\def\relbar{\mathrel{\smash-}}
\def\myline{\relbar\joinrel\relbar\joinrel\relbar\joinrel\relbar\joinrel\relbar\joinrel\relbar\joinrel\relbar
\joinrel\relbar\joinrel\relbar}
  { \ninepoint
 
\def\hiha{\vglue1pt}
 \hiha\name{P.\ Deligne} (1980), La conjecture de Weil. II. {\it Inst. Hautes {\rm \'{\it E}}tudes Sci.\
Publ.\
Math}.\   No.\ {\bf 52},  137--252.

  \hiha$\myline$
 (1982)  (Notes by J.\ S.\ Milne), Hodge cycles on abelian varieties,
in
{\it Hodge
Cycles{\rm ,} Motives{\rm ,} and Shimura Varieties\/}, pp.~9--100, {\it  Lecture Notes in
Math\/}.\ {\bf 900}, Springer-Verlag, New York.  

 \hiha$\myline$ (1989),   Le groupe fondamental de la droite projective moins trois
points, in {\it  Galois Groups over ${\scriptstyle\Bbb Q}{}$\/} (Berkeley, CA, 1987), 
79--297,
{\it Math.\ Sci.\ Res.\ Inst.\ Publ\/}.\ {\bf 16}, Springer-Verlag,
New York.

  \hiha$\myline$ 
(1990),  Cat\'{e}gories tannakiennes, in {\it The 
Grothendieck Festschrift\/}, {\it Vol\/}.\ 
II, 111--195, {\it Progr.\ Math\/}.\  {\bf 87}, 
Birkh\"{a}user Boston, Boston, MA.
\pagebreak

 \hiha\name{P.\ Deligne} and  \name{J.\ S.\ Milne} (1982),  Tannakian categories, in {\it  Hodge Cycles, 
Motives,
and Shimura Varieties\/},  pp.~101--228, {\it  Lecture Notes in
Math\/}.\ {\bf 900},
Springer-Verlag, New York.

  \hiha\name{A.\ Grothendieck} (1958),  Sur une note de Mattuck-Tate, {\it  J. Reine Angew.\
Math\/}.\ {\bf 200}, 208--215.

  \hiha$\myline$ (1969),  Standard conjectures on algebraic cycles, in {\it  
Algebraic
Geometry\/} (Internat.\ Colloq., Tata Inst.\ Fund.\ Res., Bombay, 1968),
pp.\ 
193--199, Oxford Univ. Press, London.

  \hiha\name{U.\ Jannsen} (1992),  Motives, numerical equivalence, and semi-simplicity, {\it Invent.\
Math\/}.\ {\bf 107}, 447--452.

  \hiha\name{S.\ L.\ Kleiman} (1968), Algebraic cycles and the Weil conjectures, in {\it Dix 
Expos{\rm \'{\it e}}s
sur la Cohomologie des Sch{\rm \'{\it e}}maas\/}, pp.\ 359--386, North-Holland, Amsterdam;
Masson, Paris.

  \hiha$\myline$ (1994), The standard conjectures, in {\it Motives\/}  
(Seattle, WA, 1991),
3--20, {\it Proc.\ Sympos.\ Pure Math\/}.\ {\bf 55},  Part 1, A.\ M.\ S., 
Providence,
RI.

 \hiha\name{D.\ I.\ Lieberman} (1968),  Numerical and homological equivalence of algebraic cycles
on Hodge manifolds, {\it  Amer.\ J. Math\/}.\ {\bf 90},  366--374.

  \hiha\name{J.\ S.\ Milne} (1994),  Motives over finite fields, in {\it Motives\/}  
(Seattle, WA, 1991),
401--459, {\it Proc.\ Sympos.\ Pure Math\/}.\ {\bf 55},   Part 1, A.\
M.\ S., 
Providence, RI.

  \hiha$\myline$ 
(1999a), Lefschetz classes on abelian varieties, {\it  Duke Math.\
J\/}.\ {\bf 96}, 639--675.

  \hiha$\myline$ (1999b), Lefschetz motives and the Tate conjecture, {\it  Compositio
Math\/}.\ 
{\bf 117},  45--76.

  \hiha$\myline$ (2001),  The Tate conjecture for certain abelian varieties over finite
fields, {\it  Acta Arith\/}.\ {\bf 100}, 135--166.

  \hiha\name{N.\ Saavedra Rivano} (1972),  {\it Cat{\rm \'{\it e}}gories Tannakiennes}, {\it  Lecture Notes in
Math\/}.\  {\bf 265},  Springer-Verlag, New York.

 \hiha\name{B.\ Segre} (1937), Intorno ad teorema di Hodge sulla teoria della base per le curve
di una superficie algebrica, {\it Ann.\ Mat\/}.\ {\bf 16},
157--163.

  \hiha\name{J.\ T.\ Tate} (1994), Conjectures on algebraic cycles in $l$-adic cohomology, in
{\it Motives\/}
(Seattle, WA, 1991), 71--83, {\it Proc.\ Sympos.\ Pure Math\/}.\ {\bf 55},  
Part 1, A.\ M.\ S., Providence, RI.

 \hiha\name{A.\ Weil} (1948),  {\it Vari{\rm \'{\it e}}t{\rm \'{\it e}}s Ab{\rm \'{\it e}}liennes et Courbes
Alg{\rm
\'{\it e}}briques}, {\it Actualit{\rm \'{\it e}}s Sci.\ Ind\/}.\  no.\ {\bf 1064}, Hermann \& Cie., 
Paris.

  \hiha$\myline$ (1958),  {\it Introduction {\rm \`{\it a}} l'\'{e}tude des vari{\rm \'{\it e}}t{\rm \'{\it e}}s
k{\rm \"{\it a}}hl{\rm \'{\it e}}riennes},  {\it Publications de l'Institut de Math\'{e}matique de l'Universit\'{e}
de Nancago\/}, VI. {\it Actualit{\rm \'{\it e}}s Sci.\ Ind\/}.\  no.\ {\bf 1267},  
Hermann, Paris.

 \hiha\name{T.\ Zink}
(1983),  Isogenieklassen von Punkten von Shimuramannigfaltigkeiten mit
Werten in einem endlichen K\"{o}rper, {\it  Math.\ Nachr\/}.\ {\bf 112}, 103--124.

\vglue12pt
\centerline{(Received April 3, 2001)}

 }
\end{document}